\documentclass[12pt]{article}
\usepackage{amssymb}

\newcommand{\be}{\begin{equation}}
\newcommand{\ee}{\end{equation}}
\newcommand{\bea}{\begin{eqnarray}}
\newcommand{\eea}{\end{eqnarray}}
\newcommand{\barray}{\begin{array}}
\newcommand{\earray}{\end{array}}

\newcommand{\nn}{\nonumber}
\newcommand{\bitem}{\begin{itemize}}
\newcommand{\eitem}{\end{itemize}}
\newtheorem{teo}{Theorem}[section]
\newcommand{\bt}{\begin{teo}}
\newcommand{\et}{\end{teo}}
\newtheorem{Def}{Definition}[section]
\newcommand{\bd}{\begin{Def}}
\newcommand{\ed}{\end{Def}}
\newtheorem{lem}{Lemma}[section]
\newcommand{\bl}{\begin{lem}}
\newcommand{\el}{\end{lem}}
\newtheorem{prop}{Proposition}[section]
\newcommand{\bp}{\begin{prop}}
\newcommand{\ep}{\end{prop}}
\newtheorem{cor}{Corollary}[section]
\newcommand{\bc}{\begin{cor}}
\newcommand{\ec}{\end{cor}}
\newtheorem{ex}{Example}[section]
\newcommand{\bex}{\begin{ex}}
\newcommand{\eex}{\end{ex}}
\newtheorem{rem}{Remark}[section]
\newcommand{\br}{\begin{rem}}
\newcommand{\er}{\end{rem}}

\begin{document}

\begin{center}
{\Large \textbf{Commuting ordinary differential operators \\of arbitrary genus and arbitrary rank \\ with polynomial coefficients\footnotetext[1]{This research was supported by
the Russian Foundation for Basic Research (grant nos. 11-01-00197
and 11-01-12067-ofi-m-2011) and by the ``Leading Scientific Schools'' programme (grant
no. NSh-4995-2012.1).}}}
\end{center}

\medskip
\smallskip

\begin{center}
{\large \bf {O. I. Mokhov}}
\end{center}

\medskip

\begin{center}
{\it  {To my teacher S.P.Novikov on his 75th birthday}}
\end{center}

\bigskip
\medskip

In this paper we construct examples of commuting ordinary scalar differential
operators with polynomial coefficients that are related to a spectral curve of an arbitrary genus
$g>0$ and to an arbitrary rank $r>1$ of
the vector bundle of common eigenfunctions of the commuting operators over the spectral curve. This solves completely the well-known existence problem for commuting operators of arbitrary genus and arbitrary rank with polynomial coefficients. 
The constructed commuting operators of arbitrary rank $r>1$ and arbitrary genus $g>0$ are given explicitly, they are generated by the Chebyshev polynomials $T_r (x)$.

The first operator $L_{2 r}$ from the constructed commuting pair of genus $g$ and rank $r$ is of order $2r$ and has the form
\be
L_{2r} = \left ( a T_r \left ( {d \over dx } \right ) - x^2 {d^2 \over d x^2} - 3 x {d \over dx} + x^2 + b \right )^2 - a r^2 g (g + 1) T_r \left ( {d \over dx} \right ), \label{ch}
\ee
where $T_r (x)$ is the Chebyshev polynomial of the first kind, of degree $r$, $r > 1$, the notation $T_r \left ( {d \over dx } \right )$ means the ordinary differential operator, which is the Chebyshev polynomial $T_r$ of ${d \over dx }$, $a$ is an arbitrary nonzero constant, $b$ is an arbitrary constant.
The second operator $M_{(2g + 1) r}$ from the constructed commuting pair of genus $g$ and rank $r$ is of order $(2g + 1)r$ and there exists a certain hyperelliptic relation
\be
M_{(2g + 1) r}^2 = L_{2r}^{2g + 1} + a_{2g} L_{2r}^{2g} + \ldots + a_1 L_{2r} + a_0,
\ee
where $a_i$ are some constants, $$[L_{2r}, M_{(2g + 1) r}] = 0,$$ the coefficients of the operator $M_{(2g + 1) r}$ are expressed polynomially in terms of the coefficients of the operator $L_{2r}$ and their derivatives.

To every generic point $P = (\lambda, \mu)$ of the hyperelliptic spectral curve
\be
{\mu}^2 = {\lambda}^{2g + 1} + a_{2g} {\lambda}^{2g} + \ldots + a_1 \lambda + a_0
\ee
of genus $g$ there corresponds $r$-dimensional space of common eigenfunctions $\psi (x, \lambda, \mu)$ of the constructed commuting operators $L_{2r}$ and $M_{(2g + 1) r}$ of genus $g$ and rank $r$, $$L_{2r} \psi = \lambda \psi, \ \ \ M_{(2g + 1) r} \psi = \mu \psi.$$

Recall that up to now no examples of commuting ordinary scalar differential
operators that are related to a spectral curve of genus
$g > 1$ for rank of the form $r = 6s \pm 1$, $s \geq 1$, were known. For all other values of
genus $g$ and rank $r$, explicit examples of commuting
operators even with polynomial coefficients were constructed (see [1] and references therein and here in the text below).
We conjectured in [1] that there exist commuting ordinary scalar differential
operators with polynomial coefficients that are related to a spectral curve of an arbitrary genus
$g > 1$ also for an arbitrary rank of the form $r = 6s \pm 1$, $s \geq 1$,
and in this paper we construct such examples using recent Andrey Mironov's remarkable results on self-adjoint commuting ordinary scalar differential
operators of rank 2 and an arbitrary genus $g$ (see [2], [3]). Thus,  examples of commuting ordinary scalar differential
operators with polynomial coefficients are constructed now for arbitrary genus
$g>0$ and arbitrary rank $r>1$.

The study of the commutation equation for two scalar ordinary differential operators is one of the old classical problems of the theory of ordinary differential equations (see [4]--[8]).
The problem of finding the general form of commuting operators of fixed genus $g$ and fixed rank $r$, and even constructing partial examples of such commuting operators, is very nontrivial for $g >0$ and $r > 1$. Moreover,
there was the very interesting problem of existing commuting ordinary scalar differential
operators with polynomial coefficients for fixed genus $g$ and fixed rank $r$ (in particular, this existence question was posed by I.M.Gelfand at the famous Gelfand seminar at Moscow State University in 1981). 
It is well known that commuting ordinary scalar differential
operators with polynomial coefficients give commutative subalgebras of the Weyl algebra
$W_1$, i.e., the algebra with two generators $p$ and $q$ and the relation
$[p, q] = 1$, so they are of a special interest and we also give them a special attention.
These problems are very well known and we will give here only very brief summary of the earlier obtained results.

Let us consider a system of nonlinear ordinary differential equations on the coefficients of two scalar ordinary differential operators
\be
L = \sum_{i = 0}^n u_i (x) {d^i \over d x^i}, \ \ \ M = \sum_{i = 0}^m v_i (x) {d^i \over d x^i},
\ee
that is equivalent to the commuting condition
$$[L, M] = 0.$$

The operators are usually assumed to be in the standard canonical form (Burchnall, Chaundy, [6]), i.e.
\be
u_n (x) = v_m (x) = 1, \ \ \ u_{n-1} (x) = 0.
\ee
For a pair of commuting operators this can always be achieved by a change of variables and a suitable conjugacy
(Burchnall, Chaundy, [6]).

By the Burchnall--Chaundy lemma [5], [6] any pair of commuting ordinary scalar differential
operators $L$ and $M$ is connected by a certain polynomial relation
$Q (L, M) = 0$ given by the spectral curve $Q(\lambda, \mu) = 0$ of the pair of commuting
operators:
$L \psi = \lambda \psi$, $M \psi = \mu \psi$, and common eigenfunctions of the commuting operators
define an $r$-dimensional vector bundle over the spectral curve
(the spectral curve $\Gamma$ defined by the relation $Q(\lambda, \mu) = 0$ is irreducible and is completed at infinity with a unique point $P_0$; the dimension $r$ of
the vector bundle of common eigenfunctions of a pair of commuting operators over the spectral curve
at a generic point of the curve $\Gamma$ is called {\it the rank of the
commuting pair of operators};
the rank of any pair of commuting operators is a common divisor of the orders of these
commuting operators). For commuting operators of relatively prime orders the rank is equal to 1.
In the case of rank 1 the commutation equation has been integrated in [5]--[8]. The common eigenfunctions and
the coefficients of commuting operators of rank 1 are expressed explicitly in terms of the theta-function of the spectral curve [9]. The case $l > 1$ for spectral curves of nontrivial genus $g > 0$ is much more complicated and much more interesting.

The first examples of commuting ordinary scalar differential
operators of the nontrivial ranks 2 and 3 and the nontrivial genus $g = 1$
were constructed by Dixmier [10] for the nonsingular elliptic spectral curve
$\mu^2 = \lambda^3 - \alpha$,
where $\alpha$ is an arbitrary nonzero constant:
\be
L = \left (  \left ({d \over dx}\right )^2 + x^3 + \alpha \right )^2 + 2x,
\ee
\be
M = \left (  \left ({d \over dx}\right )^2 + x^3 + \alpha \right)^3 + 3x \left ( {d \over dx} \right)^2 + 3 {d \over dx}
+ 3x (x^3 + \alpha),
\ee
where $L$ and $M$ is the commuting pair of the Dixmier operators of rank 2 and genus 1, 

$$M^2 = L^3 - \alpha, \ \ \ \ [L, M] = 0,$$
the orders of the commuting operators $L$ and $M$ are 4 and 6, rank 2; 
\be
L = \left (  \left ({d \over dx}\right )^3 + x^2 + \alpha \right )^2 + 2 {d \over dx},
\ee
\be
M = \left (  \left ({d \over dx}\right )^3 + x^2 + \alpha \right)^3 + 3 \left ( {d \over dx} \right)^4 + 3 (x^2 + \alpha) {d \over dx} + 3x,
\ee
where $L$ and $M$ is the commuting pair of the Dixmier operators of rank 3 and genus 1, 

$$M^2 = L^3 - \alpha, \ \ \ \ [L, M] = 0,$$
the orders of the commuting operators $L$ and $M$ are 6 and 9, rank 3.

These remarkable unusual examples were found by Dixmier as commutative subalgebras of the Weyl algebra $W_1$ by a quite nontrivial, purely algebraic way without any connection to the spectral theory of commuting operators [10]. Both the Dixmier examples are in the standard canonical form.

The general classification of commuting ordinary scalar differential
operators of nontrivial ranks $r > 1$ was obtained by Krichever [11]: a pair of commuting operators of rank $r$ is determined by specifying the curve $\Gamma$, a point $P_0 \in \Gamma$, a local parameter $k^{-1} (P)$ in a
neighbourhood of $P_0$, by specifying $r^2 g$ ($g$ is the genus of the curve $\Gamma$) constants $(\alpha_{ij}, \gamma_i)$, $1 \leq i \leq r g,$ $0 \leq j \leq r -2,$ called the {\it Tyurin parameters,} and by specifying $r - 1$ arbitrary functions $w_j (x)$.

The common eigenfunctions $\psi_j (x, P; x_0),$ $0 \leq j \leq r - 1,$ $P \in \Gamma,$ normalized by the condition
$$
\left. \left ({d^i \over d x^i} \psi_j (x, P; x_0) \right ) \right |_{x = x_0} = \delta_{ij}, \ \ \ 0 \leq i, j \leq r - 1,
$$
have the following analytic properties (Krichever, [11]):

1. The common eigenfunctions $\psi_j (x, P; x_0),$ $0 \leq j \leq r - 1,$ are meromorphic on the spectral curve $\Gamma$ outside $P_0$, and each one has $r g$ simple poles $\gamma_i (x_0)$, with

$$
\psi_j (x, z; x_0) \sim {\psi_{ij} (x, x_0) \over z - \gamma_i (x_0)}, \ \ \ 0 \leq j \leq r - 1, \ \ 1 \leq i \leq r g,
$$
in a neighbourhood of the pole $\gamma_i (x_0)$.

2. All the residues $\psi_{ij} (x, x_0)$ are proportional to one of them: $$\psi_{ij} (x, x_0) = \alpha_{ij} (x_0) \psi_{ir-1} (x, x_0), \ \ \ 0 \leq j \leq r - 2.$$

3. If $k^{- 1} (P)$ is a local parameter on $\Gamma$ in a neighbourhood of $P_0$, then we have the asymptotics
$$\vec{\psi} (x, P; x_0) = \left ( \sum_{s = 0}^{\infty} \vec{\xi}_s (x) k^{- s} \right ) \Phi_0 (x, k; x_0),$$
where $\vec{\psi} (x, P; x_0) = (\psi_0 (x, P; x_0), \ldots, \psi_{r - 1} (x, P; x_0))$,
$\vec{\xi}_0 (x) = (1, 0, \ldots, 0),$ $\vec{\xi}_s (x_0) = \vec{0}$ for $s \geq 1$, and $\Phi_0 (x, k; x_0) =
(\Phi_0^{ij})$ is a solution of the matrix equation $${d \Phi_0 \over d x} = S \Phi_0, \ \ \ \ \Phi^{ij}_0 (x_0, k; x_0) = \delta^{ij}, \ \ \ 0 \leq i, j \leq r - 1,$$ and
$$S = \left ( \begin{array}{cccccc} 0 & 1 & 0 & \cdots & 0 & 0 \\
0 & 0 & 1 & \cdots & 0 & 0 \\
\cdots & \cdots & \cdots & \cdots & \cdots & \cdots \\
0 & 0 & 0 & \cdots & 0 & 1 \\
k + w_0 (x) & w_1 (x) & w_2 (x) & \cdots & w_{r - 2} (x) & 0
\end{array} \right ).$$
The analytic properties 1, 2 and 3, the arbitrary constants $(\gamma_i, \alpha_{ij})$, $1 \leq i \leq r g,$ $0 \leq j \leq r -2,$ and the arbitrary functions
$w_0 (x), \ldots, w_{r - 2} (x)$ determine a vector-valued function $\vec{\psi} (x, P; x_0)$ and a commuting pair $L$, $M$ of rank $r$ and genus $g$ in general position (Krichever, [11]).

In order to find commuting operators Krichever and Novikov proposed in [11], [12] the method of deforming the Tyurin parameters $(\gamma_i (x_0), \alpha_{ij} (x_0))$ that allowed to obtain the general form of commuting ordinary scalar differential
operators of rank $r = 2$ for arbitrary elliptic spectral curve (genus $g = 1$) [13].

Let us consider the Wronskian matrix
$$ \widetilde{\psi} (x, P; x_0) = \left ( \begin{array}{cccc} \psi_0 & \psi_1 & \cdots & \psi_{r - 1} \\
\psi^{'}_0 & \psi^{'}_1 & \cdots & \psi^{'}_{r - 1} \\
\cdots & \cdots & \cdots & \cdots \\
\psi^{(r - 1)}_0 & \psi^{(r - 1)}_1 & \cdots & \psi^{(r - 1)}_{r - 1}
\end{array} \right )$$
of the vector-valued function $\vec{\psi} (x, P; x_0)$, and
$$ \widetilde{\psi}_x \widetilde{\psi}^{- 1}  = \left ( \begin{array}{ccccc}
0 & 1 & 0 & \cdots & 0 \\
0 & 0 & 1 & \cdots & 0 \\
\cdots & \cdots & \cdots & \cdots & \cdots \\
0 & 0 & 0 & \cdots & 1 \\
\chi_0 & \chi_1 & \chi_2 & \cdots & \chi_{r - 1}
\end{array} \right ),$$
where $\chi_j (x, P)$ are meromorphic functions on $\Gamma$.

For $x = x_0$ the poles of $\chi_j (x_0, P)$ coincide with $\gamma_1 (x_0), \ldots, \gamma_{rg} (x_0),$
and the ratios of the residues of the functions $\chi_j (x_0, P)$ at the points $\gamma_i (x_0)$ coincide with the parameters  $\alpha_{ij} (x_0)$: 

$$\alpha_{ij} (x_0) = {{\rm res}_{\gamma_i (x_0)} \chi_j \over {\rm res}_{\gamma_i (x_0)} \chi_{r - 1}}.$$

In a neighbourhood of $P_0$ on $\Gamma$ the functions $\chi_j (x_0, P)$ have the form
$$\chi_0 (x, P) = k + w_0 (x) + O (k^{-1}),$$
$$\chi_s (x, P) = w_s (x) + O (k^{-1}), \ \ \ 1 \leq s \leq r - 2,$$
$$\chi_{r - 1} (x, P) = O (k^{-1}).$$
The expansion of $\chi_j$ in a neighbourhood of the pole $\gamma_i (x)$ has the form
$$\chi_j (x, k) = {c_{ij} (x) \over k - \gamma_i (x)} + d_{ij} (x) + O (k - \gamma_i (x)),$$
$$c_{ij} = \alpha_{ij} c_{i r-1}, \ \ \ 0 \leq j \leq r - 1, \ 1 \leq i \leq rg.$$

Krichever and Novikov proved [11], [14] that the functional parameters $\gamma_i (x)$ and $\alpha_{ij} (x),$ $1 \leq i \leq rg,$ $0 \leq j \leq r - 2,$
satisfy the system (the Krichever--Novikov system of equations of deformation of the Tyurin parameters)
$$\gamma^{'}_i = - c_{i r-1},$$
$$\alpha^{'}_{i0} = \alpha_{i0} \alpha_{i r-2} + \alpha_{i0} d_{i r-1} - d_{i0},$$
$$\alpha^{'}_{ij} = \alpha_{ij} \alpha_{i r-2} - \alpha_{i j-1} + \alpha_{ij} d_{i r-1} - d_{ij}, \ \ \ j \geq 1.$$
For arbitrary functions $w_i (x)$, $0 \leq i \leq r - 2,$ the solution of the Krichever--Novikov system uniquely determine an array of
meromorphic functions $\chi_j (x_0, P)$ with the required analytic properties (Krichever, Novikov [11], [14]).

By the asymptotics of $\chi_j (x_0, P)$
in a neighbourhood of $P_0$ on $\Gamma$ one can reconstruct the commuting operators.
We note that the coefficients of one from a pair of commuting operators are always expressed polynomially in terms of the coefficients of the second operator of the commuting pair and their derivatives.

The general form of commuting
ordinary scalar differential operators of rank 2 for an arbitrary elliptic spectral curve was found
by Krichever and Novikov [13]. The functional parameter corresponding to the Dixmier example of rank 2, genus 1 among all
the Krichever--Novikov commuting operators of rank 2, genus 1 was found by Grinevich [15]. The class of self-adjoint
Krichever--Novikov commuting operators of rank 2, genus 1 was singled out by Grinevich and Novikov in [16]. We found explicit formulae for self-adjoint commuting operators of rank 2, genus 2 depending on two functional parameters connected by a third-order nonlinear ordinary differential relation (see [17]), but this relation was not resolved.
 The general form of commuting
ordinary scalar differential operators of rank 3 for an arbitrary elliptic spectral curve (the general commuting operators of rank 3, genus 1 are parametrized
by two arbitrary functions) was found by Mokhov [17], [18], where also the functional parameters corresponding to the Dixmier example of rank 3, genus 1 among all the commuting operators of rank 3, genus 1 were found and all commuting operators of rank 3, genus 1 with rational coefficients were singled out. Moreover, examples of commuting ordinary scalar differential
operators of genus 1 with polynomial coefficients were constructed for any rank $r$ (some commutative subalgebras of the Weyl algebra $W_1$) (see also [19] and other related papers [25]--[33]). The description of commuting ordinary scalar differential
operators with polynomial coefficients (commutative subalgebras of the Weyl algebra $W_1$) is a separate nontrivial
problem, and this problem has not been solved completely yet even for the Krichever--Novikov commuting operators of rank 2, genus 1, which are rationally parametrized by one arbitrary function (this problem was considered and studied in [20], [21]).

Recently Mironov [2] (see also earlier papers [22]--[24]) constructed for any genus $g > 1$ remarkable
examples of commuting ordinary scalar differential
operators of ranks 2 and 3 with polynomial coefficients that generalize naturally the Dixmier examples of ranks 2 and 3, genus 1:
\be
L = \left (  \left ({d \over dx}\right )^2 + x^3 + \alpha \right )^2 + g (g + 1) x,
\ee
\be
M^2 = L^{2g + 1} + a_{2g} L^{2g} + \ldots + a_1 L + a_0,
\ee
where $a_i$ are some constants, $\alpha$ is an arbitrary nonzero constant, $L$ and $M$ are the Mironov commuting operators of rank 2, genus $g$ (the orders of the operators $L$ and $M$ are 4 and $4g + 2$, respectively), the coefficients of the operator $M$ are expressed polynomially in terms of the coefficients of the operator $L$ and their derivatives, $[L, M] = 0$;
\be
L = \left (  \left ({d \over dx}\right )^3 + x^2 + \alpha \right )^2 + g (g + 1) {d \over dx},
\ee
\be
M^2 = L^{2g + 1} + a_{2g} L^{2g} + \ldots + a_1 L + a_0,
\ee
where $a_i$ are some constants, $\alpha$ is an arbitrary nonzero constant, $L$ and $M$ are the Mironov commuting operators of rank 3, genus $g$ (the orders of the operators $L$ and $M$ are 6 and $6g + 3$, respectively), the coefficients of the operator $M$ are expressed polynomially in terms of the coefficients of the operator $L$ and their derivatives, $[L, M] = 0$.

Using Mironov's results, in our paper [1] we constructed examples of commuting ordinary scalar differential
operators with polynomial coefficients that are related to a spectral curve of an arbitrary genus
$g$ for an arbitrary even rank $r = 2k$, $k > 1$, and for an arbitrary rank of the form $r = 3k$, $k \geq 1$.

The operators $L$ and $M$ of orders $4k$ and $4kg + 2k$, respectively,
\be
L = \left (  \left ({d \over dx}\right )^{2k} - 2x \left ({d \over dx}\right )^{k} - k \left ({d \over dx}\right )^{k - 1} + \left ({d \over dx}\right )^3 + x^2 + \alpha \right )^2 + g (g + 1) {d \over dx},
\ee
\be
M^2 = L^{2g + 1} + a_{2g} L^{2g} + \ldots + a_1 L + a_0,
\ee
where $a_i$ are some constants, $\alpha$ is an arbitrary nonzero constant, are commuting operators of rank $r = 2k$, $k > 1$, genus $g$, $[L, M] = 0$, the coefficients of the operator $M$ are expressed polynomially in terms of the coefficients of the operator $L$ and their derivatives [1]. For $k >2$ the commuting operators $L$ and $M$ have the standard canonical form.

The operators $L$ and $M$ of orders $6k$ and $6kg + 3k$, respectively,
\bea
&& L = \left (  \left ({d \over dx}\right )^{3k} - 3x \left ({d \over dx}\right )^{2k} - 3k \left ({d \over dx}\right )^{2k - 1} + 3 x^2 \left ({d \over dx}\right )^k + 3 k x \left ({d \over dx}\right )^{k - 1} + \right. \nn\\ && \left. + k (k - 1) \left ({d \over dx}\right )^{k - 2} + \left ({d \over dx}\right )^2 - x^3 + \alpha \right )^2 - g (g + 1) x,
\eea
\be
M^2 = L^{2g + 1} + a_{2g} L^{2g} + \ldots + a_1 L + a_0,
\ee
where $a_i$ are some constants, $\alpha$ is an arbitrary nonzero constant, are commuting operators of rank $3k$, $k \geq 1$, genus $g$, $[L, M] = 0$, the coefficients of the operator $M$ are expressed polynomially in terms of the coefficients of the operator $L$ and their derivatives [1]. For $k >1$ the commuting operators $L$ and $M$ have the standard canonical form.

In [2] Mironov proved that for commuting operators $L,$ $M$ of rank 2 with hyperelliptic spectral curve $\Gamma$ of arbitrary genus $g$,
\be
\Gamma: \ {\mu}^2 = {\lambda}^{2g + 1} + a_{2g} {\lambda}^{2g} + \ldots + a_1 {\lambda} + a_0, \label{gamma}
\ee
where $a_i$ are some constants, 
the operator $L$ of order 4 is self-adjoint if and only if

$$\chi_1 (x, P) = \chi_1 (x, \sigma (P)),$$
where $\sigma$ is the involution
$$\sigma (\lambda, \mu) = (\lambda, - \mu)$$
on the hyperelliptic curve $\Gamma$.

Any self-adjoint operator $L$ of order 4 in the standard canonical form can be represented as

\be
L = \left ( {d^2 \over d x^2} + V (x)\right)^2 + W (x). \label{op}
\ee

By the Mironov theorem [2], if the operator $L$ of order 4 from a pair of commuting operators $L,$ $M$ of rank 2 with hyperelliptic spectral curve $\Gamma$ (\ref{gamma}) of arbitrary genus $g$ is self-adjoint, then
the operator $L$, obviously, has the form (\ref{op}) and the corresponding functions $\chi_0 (x, P),$ $\chi_1 (x, P)$ have the form
\be
\chi_0 (x, P) = - {1 \over 2} {Q_{xx} \over Q (x, \lambda)} + {\mu \over Q (x, \lambda)} - V (x), \ \ \ \
\chi_1 (x, P) = {Q_x \over Q (x, \lambda)},
\ee 
where $P = (\lambda, \mu)$ and $Q (x, \lambda)$ is a polynomial in $\lambda$ of degree $g$ with coefficients depending on $x$,
\be
Q (x, \lambda) = (\lambda - \lambda_1 (x)) \cdots (\lambda - \lambda_g (x)),
\ee
such that the following relation holds:
\bea
&& (\lambda - W(x)) (Q (x, \lambda))^2 - V(x) (Q_{x})^2 + {1 \over 4} (Q_{xx})^2 - {1 \over 2} Q_{x} Q_{xxx}
+ \nn\\ && + Q \left (V_{x} Q_{x} + 2 V(x) Q_{xx} + {1 \over 2} Q_{xxxx} \right ) =
{\lambda}^{2g + 1} + a_{2g} {\lambda}^{2g} + \ldots + a_1 {\lambda} + a_0.
\eea

These results allowed to prove the following theorem (Mironov, [3]):
the self-adjoint operator 
\be
L = \left ( {d^2 \over dx^2} + \alpha_1 S (x) + \alpha_0 \right )^2 +
\alpha_1 c_2 g (g + 1) S (x), 
\ee
where the function $S(x)$ satisfies the equation
\be
(S_{x})^2 = c_2 (S(x))^2 + c_1 S(x) + c_0,
\ee
$\alpha_1$ and $c_2$ are arbitrary nonzero constants, $\alpha_0$, $c_1$ and $c_0$
are arbitrary constants, form a pair of commuting operators $L,$ $M$ of rank 2 with a certain hyperelliptic spectral curve $\Gamma$ (\ref{gamma}) of genus $g$.

In particular, if $S(x) = \cosh x$, $c_2 = 1$, $c_1 = 0,$ $c_0 = - 1$, then all the conditions are satisfied and
the operator
\be
L = \left ( {d^2 \over dx^2} + \alpha_1 \cosh (x) + \alpha_0 \right )^2 +
\alpha_1 g (g + 1) \cosh x, \label{l}
\ee
is the first operator of order 4 from the commuting pair $L$, $M$ of genus $g$ and rank 2 with a hyperelliptic spectral curve.

Let us consider the change of variable
$$x = \ln P (z),$$
where $P (z)$ is an arbitrary nonconstant function.
Then the operator (\ref{l}) is transformed to the form
\bea
&&
L = \left ( {(P (z))^2 \over (P_{z})^2} {d^2 \over dz^2} + {P (z) ((P_{z})^2 - P (z) P_{zz}) \over (P_{z})^3} {d \over dz} + \alpha_1 {(P (z))^2 + 1 \over 2 P (z)} + \alpha_0 \right )^2 + \nn\\ && +
\alpha_1 g (g + 1) {(P (z))^2 + 1 \over 2 P (z)}. \label{l2}
\eea
In particular, if $P (z)$ is an arbitrary nonconstant rational function, then we get an operator of order 4 from the commuting pair of genus $g$ and rank 2 with rational coefficients, but this operator is not in the standard canonical form.

Let us consider the special case
$$P (z) = (z + \sqrt {z^2 - 1})^r, \ \ \ \ r = \pm 1, \pm 2, \ldots $$

Then the operator (\ref{l2}) can be represented as
\be
L = \left ( (1 - z^2){d^2 \over dz^2} - z {d \over dz} + a T_r (z) + b \right )^2 -
a r^2 g (g + 1) T_r (z), \label{l3}
\ee
where $r$ is an arbitrary nonvanishing integer, 
$T_r (z)$ is the Chebyshev polynomial of the first kind, of degree $r$ or $-r$, $a$ is an arbitrary nonzero constant,
$b$ is an arbitrary constant.
The operator (\ref{l3}) is the first operator of order 4 from the commuting pair $L$, $M$ of genus $g$ and rank 2 with a hyperelliptic spectral curve for arbitrary nonvanishing integer $r$, but it is not in the standard canonical form. We note that all these operators have polynomial coefficients and give commutative subalgebras of the Weyl algebra $W_1$.

Recall that 
\be
T_0 (z) = 1, \ \ \ T_1 (z) = z, \ \ \ T_r (z) = 2 z T_{r-1} (z) - T_{r-2} (z), \ \ \ T_{- r} (z) = T_{r} (z). 
\ee
It is very interesting that the Chebyshev polynomials of the first kind $T_{r} (z)$ are commuting polynomials
\be
T_n (T_m (z)) = T_{nm} (z) = T_m (T_n (z)).
\ee
We consider that it is no mere chance. Earlier we noted an important role the Chebyshev polynomials of the first kind in the theory of commuting operators with polynomial coefficients and with elliptic spectral curve (see [21]).
It would be very interesting to clarify these interconnections. 

After a natural automorphism of the Weyl algebra $W_1$ we obtain the operator (\ref{ch}) from a commuting pair of rank $r$ and genus $g$ with polynomial coefficients in the standard canonical form.

{\bf Theorem.}  
{\it The operators $L_{2r}$ and $M_{(2g + 1) r}$ of orders $2r$ and $(2g + 1) r$, respectively,
\be
L_{2r} = \left ( a T_r \left ( {d \over dx } \right ) - x^2 {d^2 \over d x^2} - 3 x {d \over dx} + x^2 + b \right )^2 - a r^2 g (g + 1) T_r \left ( {d \over dx} \right ), \label{ch2}
\ee
\be
M_{(2g + 1) r}^2 = L_{2r}^{2g + 1} + a_{2g} L_{2r}^{2g} + \ldots + a_1 L_{2r} + a_0,
\ee
where $a_i$ are some constants, $T_r (x)$ is the Chebyshev polynomial of the first kind, of degree $r$, $r > 1$, the notation $T_r \left ( {d \over dx } \right )$ means the ordinary differential operator, which is the Chebyshev polynomial $T_r$ of ${d \over dx }$, $a$ is an arbitrary nonzero constant, $b$ is an arbitrary constant, are commuting operators of rank $r$, genus $g$, $[L_{2r}, M_{(2g + 1) r}] = 0$, the coefficients of the operator $M_{(2g + 1) r}$ are expressed polynomially in terms of the coefficients of the operator $L_{2r}$ and their derivatives. For $r>3$ the commuting operators $L_{2r}$ and $M_{(2g + 1) r}$ have the standard canonical form {\rm (}for $a = 1/2^{r-1}${\rm )}. For $r = 1$ this pair of operators is commuting one of rank 2 and genus $g$.}

{\bf Examples.}

1) Rank 4, genus $g$:

\bea
 L_8 =&&\left (  \left ({d \over dx}\right )^4 - (x^2 + 1)  \left ({d \over dx}\right )^2 - 3x  {d \over dx} + x^2 + \alpha \right )^2 - \nn\\ && - 16 g (g + 1) \left (  \left ({d \over dx}\right )^4 -  \left ({d \over dx}\right )^2\right ).
\eea

2) Rank 5, genus $g$:

\bea
 L_{10} =&&\left (  \left ({d \over dx}\right )^5 - {5 \over 4} \left ({d \over dx}\right )^3 - x^2  \left ({d \over dx}\right )^2 - \left (3x - {5 \over 16} \right ) {d \over dx} + x^2 + \alpha \right )^2 - \nn\\ && - 25 g (g + 1) \left (  \left ({d \over dx}\right )^5 - {5 \over 4} \left ({d \over dx}\right )^3 + {5 \over 16} {d \over dx}\right ).
\eea

3) Rank 6, genus $g$:

\bea
 L_{12} =&&\left (  \left ({d \over dx}\right )^6 - {3 \over 2} \left ({d \over dx}\right )^4 - \left (x^2 - {9 \over 16} \right )  \left ({d \over dx}\right )^2 - 3x {d \over dx} + x^2 + \alpha \right )^2 - \nn\\ && - 36 g (g + 1) \left (  \left ({d \over dx}\right )^6 - {3 \over 2} \left ({d \over dx}\right )^4 + {9 \over 16}  \left ( {d \over dx}\right )^2 \right ).
\eea

4) Rank 7, genus $g$:

\bea
 L_{14} =&&\left (  \left ({d \over dx}\right )^7 - {7 \over 4} \left ({d \over dx}\right )^5 + {7 \over 8} \left ({d \over dx}\right )^3 - x^2  \left ({d \over dx}\right )^2 - \left (3x + {7 \over 64} \right ) {d \over dx} + x^2 + \alpha \right )^2 - \nn\\ && - 49 g (g + 1) \left (  \left ({d \over dx}\right )^7 - {7 \over 4} \left ({d \over dx}\right )^5 + {7 \over 8} \left ({d \over dx}\right )^3 -
  {7 \over 64} {d \over dx}\right ).
\eea

\smallskip

\begin{center}
\bf {References}
\end{center}

\smallskip

[1] O.I. Mokhov. On commutative subalgebras of the Weyl algebra that are related to commuting operators of arbitrary rank and genus. arXiv:1201.5979.

[2] A.E. Mironov. Self-adjoint commuting differential operators and commutative subalgebras of the Weyl algebra.
arXiv: 1107.3356. 

[3] A.E. Mironov. Periodic and rapid decay rank two self-adjoint commuting differential operators. arXiv:1302.5735. 

[4] E.L. Ince. Ordinary differential equations. Longman, Green \& Co., London, 1926; reprint: Dover, New York, 1944.

[5] J.L. Burchnall, I.W. Chaundy. Commutative ordinary differential operators. {\it Proc. London Math. Society}, ser. 2, {\bf 21} (1923), 420--440.

[6] J.L. Burchnall, I.W . Chaundy.  Commutative ordinary differential operators. {\it Proc. Royal Soc. London}, ser. A, {\bf 118} (1928), 557--583.

[7] J.L. Burchnall, I.W. Chaundy.  Commutative ordinary differential operators. {\it Proc. Royal Soc. London}, ser. A, {\bf 134} (1931), 471--485.

[8] H.F. Baker. Note on [6]. {\it Proc. Royal Soc. London}, ser. A, {\bf 118} (1928), 584--593.

[9] I.M. Krichever. Integration of nonlinear equations by the methods of algebraic geometry. {\it Funktsional. Anal. i Prilozhen.}, {\bf 11}:1 (1977), 15--31 (in Russian); English translation in {\it Functional Anal. Appl.}, {\bf 11}:1 (1977), 12--26.

[10] J. Dixmier. Sur les alg\`ebres de Weyl. {\it Bull. Soc. Math. France}, {\bf 96} (1968), 209--242.

[11] I.M. Krichever. Commutative rings of ordinary linear differential operators. {\it Funktsional. Anal. i Prilozhen.}, {\bf 12}:3 (1978), 20--31 (in Russian); English translation in {\it Functional Anal. Appl.}, {\bf 12}:3 (1978), 175--185.

[12] I.M. Krichever, S.P. Novikov. Holomorphic bundles over Riemann surfaces and the Kadomtsev--Petviashvili equation. I. {\it Funktsional. Anal. i Prilozhen.}, {\bf 12}:4 (1978), 41--52 (in Russian); English translation in {\it Functional Anal. Appl.}, {\bf 12}:4 (1978), 276--286.

[13] I.M. Krichever, S.P. Novikov. Holomorphic bundles and nonlinear equations. Finite-zone solutions of rank 2. {\it Dokl. Akad. Nauk SSSR}, {\bf 247}:1 (1979), 33--36 (in Russian); English translation in {\it Soviet Math. Dokl.}.

[14] I.M. Krichever, S.P. Novikov. Holomorphic bundles over algebraic curves and non-linear equations. {\it Uspekhi Matem. Nauk}, {\bf 35}:6 (1980), 47--68 (in Russian); English translation in {\it Russian Math. Surveys}, {\bf 35}:6 (1980), 53--79.

[15] P.G. Grinevich. Rational solutions for the equation of commutation of differential operators. {\it Funktsional. Anal. i Prilozhen.}, {\bf 16}:1 (1982), 19--24 (in Russian); English translation in {\it Functional Anal. Appl.}, {\bf 16}:1 (1982), 15--19.

[16]  S.P. Novikov, P.G. Grinevich. Spectral theory of commuting operators of rank two with periodic coefficients. {\it Funktsional. Anal. i Prilozhen.}, {\bf 16}:1 (1982), 25--26 (in Russian); English translation in {\it Functional Anal. Appl.}, {\bf 16}:1 (1982), 19--20.

[17] O.I. Mokhov. Commuting differential operators of rank 3, and nonlinear differential equations. {\it Izvestiya AN SSSR, ser. matem.}, {\bf 53}:6 (1989), 1291--1315 (in Russian); English translation in {\it Math. USSR, Izvestiya}, {\bf 35}:3 (1990), 629--655.

[18] O.I. Mokhov. Commuting ordinary differential operators of rank 3 corresponding to an elliptic curve. {\it Uspekhi Matem. Nauk}, {\bf 37}:4 (1982), 169--170 (in Russian); English translation in {\it Russian Math. Surveys}, {\bf 37}:4 (1982), 129--130.

[19] P. Dehornoy. Op\'erateurs diff\'erentiels et courbes elliptiques. {\it Compositio Math.}, {\bf 43}:1 (1981), 71--99.

[20] O.I. Mokhov. On commutative subalgebras of Weyl algebra, which are associated with an elliptic curve.
International Conference on Algebra in Memory of A.I.Shirshov (1921--1981). Barnaul, USSR, 20--25 August 1991.
Reports on theory of rings, algebras and modules. 1991. P. 85.

[21] O.I. Mokhov. On the commutative subalgebras of Weyl algebra, which are generated by the Chebyshev polynomials.
Third International Conference on Algebra in Memory of M.I.Kargapolov (1928--1976). Krasnoyarsk, Russia, 23--28 August 1993. Krasnoyarsk: Inoprof, 1993. P. 421.

[22]  A.E. Mironov. On commuting differential operators of rank 2. {\it Sib. Elektron. Matem. Izvestiya [Siberian Electronic Mathematical Reports]}, {\bf 6} (2009), 533--536 (in Russian).

[23] A.E. Mironov. Commuting rank 2 differential operators corresponding to a curve of genus 2. {\it Funktsional. Anal. i Prilozhen.}, {\bf 39}:3 (2005), 91--94 (in Russian); English translation in
{\it Functional Anal. Appl.}, {\bf 39}:3 (2005), 240--243.

[24] A.E. Mironov. A ring of commuting differential operators of rank 2 corresponding to a curve of genus 2.
{\it Matem. Sbornik}, {\bf 195}:5 (2004), 103--114 (in Russian); English translation in
{\it Sbornik: Mathematics}, {\bf 195}:5 (2004), 711--722.

[25] A.E. Mironov. Commuting higher rank ordinary differential operators. 

\noindent
arXiv:1204.2092.

[26] O.I. Mokhov. Canonical Hamiltonian representation of the Krichever-Novikov equation. {\it Matem. Zametki}, {\bf 50}:3 (1991), 939--945 (in Russian); English translation in {\it Mathematical Notes}, {\bf 50}:3 (1991), 939--945.

[27] F. Grunbaum, Commuting pairs of linear ordinary differential operators of orders four and
six. {\it Phys. D}, {\bf 31}:3 (1988), 424--433.

[28] D. Mumford. An algebro-geometric constructions of commuting operators and of solutions to
the Toda lattice equations, Korteweg--de Vries equations and related non-linear equations. In
Proc. Internat. Symp. on Alg. Geom., Kyoto 1977, Kinokuniya Publ. (1978) 115--153.

[29] E. Previato. Seventy years of spectral curves: 1923--1993. {\it Integrable Systems and Quantum Groups.
Lecture Notes in Mathematics}, {\bf 1620}, 1996, 419-481.

[30] F. Gesztesy, R. Weikard. A characterization of all elliptic algebro-geometric solutions of the AKNS hierarchy.
{\it Acta Mathematica}, {\bf 181}:1 (1998), 63--108.

[31] G.A. Latham. The computational approach to commuting ordinary differential operators of orders six and nine.
{\it The Journal of the Australian Mathematical Society. Series B. Applied Mathematics}, {\bf 35}:04 (1994), 399--419.

[32] G.A. Latham, E. Previato. Darboux transformations for higher-rank Kadomtsev--Petviashvili and Krichever--Novikov equations. {\it Acta Applicandae Mathematica}, {\bf 39}:1-3 (1995), 405--433.

[33] D. Zuo. Commuting Differential Operators of Rank 3 Associated to a Curve of Genus 2.
{\it SIGMA}, {\bf 8}:044 (2012), 11 pp. arXiv:1105.5774. 

\begin{flushleft}
%{\bf O.I. Mokhov}\\
Department of Geometry and Topology,\\
Faculty of Mechanics and Mathematics,\\
Lomonosov Moscow State University\\
Moscow, 119991 Russia\\
{\it E-mail\,}: mokhov@mi.ras.ru; mokhov@landau.ac.ru; mokhov@bk.ru\\
\end{flushleft}

\bigskip
\bigskip
\medskip
\smallskip

\begin{center}
{\large \bf {Abstract}}
\end{center}

\bigskip
\bigskip

In this paper we construct examples of commuting ordinary scalar differential
operators with polynomial coefficients that are related to a spectral curve of an arbitrary genus
$g>0$ and to an arbitrary rank $r>1$ of
the vector bundle of common eigenfunctions of the commuting operators over the spectral curve. This solves completely the well-known existence problem for commuting operators of arbitrary genus and arbitrary rank with polynomial coefficients.
The constructed commuting operators of arbitrary rank $r>1$ and arbitrary genus $g>0$ are given explicitly, they are generated by the Chebyshev polynomials $T_r (x)$.

\end{document}